\documentclass[12pt, a4article]{amsart}   	
\usepackage[top=1.17in, bottom=1.17in, left=1.19in, right=1.19in]{geometry}

\usepackage{amssymb,amscd,graphics, amsfonts, multicol}
\usepackage{graphics}
\usepackage{stmaryrd}
\usepackage{DotArrow}
\usepackage{amsmath, amsthm,wasysym}
\usepackage{epsf}
\usepackage{xypic}
\usepackage{tikz-cd}

\def\0{{\theta}}

\DeclareMathOperator{\supp}{supp}

\DeclareMathOperator{\Q}{\mathbb{Q}}
\DeclareMathOperator{\LS}{LS}
\DeclareMathOperator{\wt}{wt}
\DeclareMathOperator{\ide}{{\it id}}

\title[Seshadri stratifications and Schubert varieties]{Seshadri stratifications and Schubert varieties: {a geometric construction of a standard monomial theory}}
\author{Rocco Chiriv\`i}
\address{Dipartimento di Matematica e Fisica ``Ennio De Giorgi'', Universit\`a del Salento, Lecce, Italy}
\email{rocco.chirivi@unisalento.it}

\author{Xin Fang} 
\address{Department Mathematik/Informatik, Universit\"at zu K\"oln, 50931, Cologne, Germany}
\email{xinfang.math@gmail.com}

\author{Peter Littelmann}
\address{Department Mathematik/Informatik, Universit\"at zu K\"oln, 50931, Cologne, Germany}
\email{peter.littelmann@math.uni-koeln.de}

\theoremstyle{plain}
\newtheorem{theorem}{Theorem}[section]
\newtheorem{coro}[theorem]{Corollary}
\newtheorem{proposition}[theorem]{Proposition}
\newtheorem{lemma}[theorem]{Lemma}

\theoremstyle{definition}

\newtheorem{remark}[theorem]{Remark}
\newtheorem{conclusion}[theorem]{Conclusion}

\newtheorem{definition}[theorem]{Definition}
\makeindex



\begin{document}

\begin{abstract}
A standard monomial theory for Schubert varieties is constructed exploiting (1) the geometry of the Seshadri stratifications of Schubert varieties by their Schubert subvarieties and (2) the combinatorial LS-path character formula for Demazure modules. The general theory of Seshadri stratifications is improved by using arbitrary linearization of the partial order and by weakening the definition of balanced stratification.
\end{abstract}

\maketitle
\begin{center}
\textit{A Corrado, amico e maestro}
\end{center}

\section{Introduction}

\subsection{Standard monomial theory.} There is no formal definition of what a standard monomial theory is. If $\mathbb{K}$ is a field and $A$ a $\mathbb{K}$--algebra, one tries to locate ring generators for $A$ and a class of monomials in these generators, called standard monomials, which forms a basis of $A$ as a $\mathbb{K}$--vector space. The relations expressing the \emph{non}-standard monomials as linear combinations of standard ones are called the straightening relations of the theory. Usually the indexing set of the generators are combinatorial in nature, typically in terms of a partially ordered set. Also the rules to decide whether a monomial is standard or not is combinatorial in terms of the same partially ordered set.

One of the first examples of a standard monomial theory is the study of Hodge of the Grassmannians and their Schubert varieties (see \cite{Ho1}, \cite{Ho2}). In modern terms, the Pl\"ucker coordinates and the Pl\"ucker relations define a structure of a Hodge algebra on the coordinate ring of the cone over the Pl\"ucker embedding of a grassmannian. This structure, algebraically defined in \cite{DEP} by De Concini, Eisenbud and Procesi, has some geometric consequences: the projective normality and the Cohen-Macaulay property for the Schubert varieties in the grassmannians and for determinantal varieties (see \cite{Igusa}, \cite{Musili}).

However, the notion of Hodge algebra was too rigid to treat, at least in a natural way, the case of an arbitrary Schubert variety in a flag variety $G/P$ with $G$ a semisimple, simply connected algebraic group and $P$ a parabolic subgroup. One possible generalization was enabled by: (a) the introduction in \cite{L2} (see also \cite{L3}) of the combinatorial tool of Lakshmibai-Seshadri paths, LS-paths for short, for a $G$--module $V(\lambda)$, with $\lambda$ a dominant weight of $G$ with stabilizer $P$ and (b) the construction of a suitable basis $\{x_\pi\,|\,\pi\in\LS^+_1(\lambda)\}$ of $V(\lambda)^*$ in \cite{L1} (see also \cite{LLM}) parametrized by the set $\LS^+_1(\lambda)$ of LS-paths of shape $\lambda$.

This generalization, termed LS-algebra (see \cite{Ch}), has generators indexed by LS-paths over a partially ordered set, a notion of standardness for monomials in these generators and straightening relations as a Hodge algebra. The general theory of LS-algebras gives a degeneration result and allows a new proof, in the spirit of the standard monomial theory, of the projectively normality and of the Cohen-Macaulay property for Schubert varieties in $G/P$. The basis $x_\pi$, $\pi\in\LS^+_1(\lambda)$, have been fruitfully used to construct standard monomial theories for other class of varieties with slightly modified version of LS-algebras (see \cite{Br_wonderful, Ch_multicone, Ch_ring_CSV, Ch_equations_SV, Ch_model}).

Of course, on the way from from Hodge algebras to LS-algebras there is the fundamental work of Lakshmibai, Musili and Seshadri on standard monomial theory, leading to Lakshmibai's conjecture \cite{LS2} on what is now termed the LS-path character formula. A very good account on the  case of Grassmannians, generalizations to other semisimple algebraic groups and applications can be found 
in the  book \cite{SMT} by Seshadri. Historical accounts and overviews are given in \cite{LS2,S2, S3}. One of the aims of our approach towards a standard monomial theory for Schubert varieties is to give a construction which is uniform for all types, purely geometric (up to a character formula), but where the standard monomials have the same desirable properties (for example compatibility with all Schubert varieties) as in the approach by Lakshmibai, Musili and Seshadri.

\subsection{The search of a geometric 
interpretation: the case of Grassmannians.}
Though the notion of an LS-path is purely combinatorial in nature, it has a geometric background (see \cite{LSII,LSIII,LSIV,LS}). 
And the definition of the basis formed by the functions $x_\pi$, $\pi\in\LS^+_1(\lambda)$, uses in an essential way the theory of quantum groups 
at a root of unity. It seems quite natural to search for a geometric meaning and construction of both: 
the LS-paths and the path vectors. 
As we will see below, a possible route to such a geometric interpretation is a suitable generalization of the theory of Newton-Okounkov body: instead of a \emph{flag} we consider a \emph{net} of subvarieties. This project has been started in \cite{FL} by the second and third author where they studied the case of Grassmannians. Already this case is non-trivial and gave a very interesting starting point; let us briefly recall it.

Let $R$ be the coordinate ring of the cone over the Grassmannian $G_{d,n}$ of $d$--dimensional vector spaces in $\mathbb{C}^n$ embedded in a projective space by the Pl\"ucker coordinates. For each maximal chain $\mathfrak{C}$ of Schubert varieties $X_{\tau_0}\subsetneq X_{\tau_{1}}\subsetneq\cdots\subsetneq X_{\tau_r}$ in $G_{d,n}$ one can define a valuation $\mathcal{V}_\mathfrak{C}$ on $R$ with values in $\mathbb{N}^{r+1}$. This valuation has at most one-dimensional leaves and the Newton-Okounkov body of $R$ associated to $\mathcal{V}_\mathfrak{C}$ is unimodular to the Gelfand-Tsetlin polytope.

In order to have a ``global'' object, one can define a quasi-valuation $\mathcal{V}$ by taking the minimum over all maximal chains of the valuations $\mathcal{V}_\mathfrak{C}$ with respect to a linearization of the partial order given by the containment of Schubert varieties. The quasi-valuation $\mathcal{V}$ has the following properties: (i) it has leaves that are at most one-dimensional; (ii) the Pl\"ucker coordinates, i.e. the generators $R$ as a Hodge algebra, are representatives for the non-zero leaves of degree $1$; (iii) the quasi-valuation is additive on a monomial if and only if the monomial is standard; (iv) the straightening relations of $R$ are just expressions of non-standard monomials as sums of ``greater'' leaves. In particular the Hodge algebra structure on $R$ is defined in terms of certain valuations and a quasi-valuation.

\subsection{LS-algebras and valuations} 
The second step of our project (see \cite{CFL2}) showed that an LS-algebra may be defined in terms of valuations and a quasi-valuation as well. Also in this case one recovers the generators of the algebra via valuations having at most one-dimensional leaves, the quasi-valuation is additive on standard monomials and the straightening relations have a clear interpretation in terms of the quasi-valuation. However this approach is suitable only for a variety whose homogeneous coordinate ring admits an LS-algebra structure; moreover, in the case of Schubert varieties it uses \emph{a priori} the combinatorial LS-paths and the basis constructed via quantum groups (\cite{L1}).

\subsection{Seshadri stratifications} The next step was the introduction of the notion of a Seshadri stratification in \cite{CFL}. A Seshadri stratification for the embedded projective variety $X\subseteq\mathbb{P}(V)$ is the data of a collection of subvarieties $X_p$ and regular functions $f_p\in\mathbb{K}[\hat{X}]$ on the cone $\hat{X}$ of $X$ in $V$, with $p$ in a finite partially ordered set $A$ having a unique maximal element $p_{\max}$ with $X_{p_{\max}} = X$. These subvarieties and functions must satisfy: (1) each subvariety $X_p$ is smooth in codimension $1$ and for $p>q$ a covering relations in $A$, the variety $X_q$ is of codimension $1$ in $X_p$; (2) for any $q\not\leq p$ the function $f_q$ vanishes on $X_p$; (3) $f_p$ vanishes in $X_p$ exactly on the (set theoretic) union of all divisors $X_q$, $q\in A$.

This general setting has a distinctive character which did not show up in the Grassmannian case nor in the theory of Newton-Okounkov body. A function $f_p$ is not necessary a uniformizer in the local ring $\mathcal{O}_{\hat{X}_p,\hat{X}_q}$, with $p>q$ a covering in $A$; indeed taking $f_p$ as a uniformizer is not always possible in our context since, in general, $p$ may cover many elements in $A$. So we have the notion of the bond $b_{p,q}$, a natural number giving the order of vanishing of $f_p$ along its divisor $X_q$. Taking into account such multiplicities, we construct a valuation $\mathcal{V}_\mathfrak{C}$ on $R$ having values in $\Q^{\mathfrak{C}}\subseteq\Q^A$ (no longer in $\mathbb{Z}^\mathfrak{C}$) with a procedure similar to that of the Newton-Okounkov theory and influenced by the idea of Allen Knutson \cite{knutson} to use Rees valuations.

As above, one can use this valuation to construct a quasi-valuation $\mathcal{V}$ on $R$ with values in $\Q^A$ by taking the minimum with respect the lexicographic order on $\Q^A$ induced by a fixed linearization $\leq^t$ of the partial order of $A$. This quasi-valuation has leaves that are at most one-dimensional and it takes non-negative values, i.e. its image $\Gamma$ is contained in $\Q^A_{\geq0}$. In particular $\Gamma$ is the union of finitely generated monoids $\Gamma_\mathfrak{C}$, one for each maximal chain $\mathfrak{C}$ in $A$; this is called a fan of monoids. One can define a Newton-Okounkov simplicial complex and an integral structure associated to $\mathcal{V}$ having  properties similar to those of the usual Newton-Okounkov body of a valuation. Further $\mathcal{V}$ induces a flat degeneration of $X$ to a reduced union $X_0$ of projective (not necessarily normal) toric varieties. Moreover, $X_0$ is equidimensional, and its irreducible components are in bijection with the maximal chains in $A$.

One can use the quasi-valuation $\mathcal{V}$ to set up a standard monomial theory for the homogeneous coordinate ring $R$ of $X$. Such a theory is very similar to the one of an LS-algebra in case the stratification is normal and balanced. When all the monoids $\Gamma_\mathfrak{C}$ are saturated, i.e. $\mathbb{Z}\Gamma_\mathfrak{C}
\cap\Q^\mathfrak{C}_{\geq 0} = \Gamma_\mathfrak{C}$, we say that the Seshadri stratification is normal. If $\Gamma$ is independent of the linearization $\leq^t$, where $\leq^t$ runs in a sufficiently large class of linearizations, then the stratification is called balanced.

When the stratification is normal each irreducible component of the semi-toric variety $X_0$ is a normal toric variety. In our paper \cite{CFL4} we study some properties of normal Seshadri stratifications. In particular, we show that a Gr\"obner basis of the defining ideal of $X_0$ can be lifted to define the embedded projective variety and we discuss the Koszul and Gorenstein properties for $R$.

\subsection{Seshadri stratifications for Schubert varieties.}
Let us turn back to the context of Schubert varieties in $G/P$. Let $W$ be the Weyl group of $G$, let $W_P$ be the Weyl group of the Levi subgroup associated to $P$ and fix also a dominant weight $\lambda$ with stabilizer $P$. The Schubert varieties in $G/P$ are indexed by $W/W_P$; this set is partially ordered by the Bruhat order. Now fix a Schubert variety $X(\tau)$, $\tau\in W/W_P$, and consider its embedding in $\mathbb{P}(V(\lambda)_\tau)$ where $V(\lambda)_\tau\subseteq V(\lambda)$ is the Demazure module of $\tau$. In \cite{CFL3} we have proved that the collection of Schubert subvarieties contained in $X(\tau)$ and of the functions of extremal weight in $V(\lambda)_\tau^*$ defines a Seshadri stratification for $X(\tau)$. Moreover, the main point in \cite{CFL3} is to show that the standard monomial theory defined in \cite{L1} fits into the concept of our theory. More precisely, let $\mathcal{V}$ be the quasi-valuation arising from this Seshadri stratification. 
We show that the set of LS-paths $\LS^+(\lambda)$ coincides with the fan of monoids $\Gamma$ of $\mathcal{V}$ and thus give an algebro-geometric interpretation of the LS-path model of a representation. 
In addition we show that the standard monomials defined in \cite{L1} are indeed representatives of the non-zero 
leaves of the quasi-valuation $\mathcal{V}$ for this Seshadri stratification, and the notion of the \textit{standardness} of a product in \cite{CFL3} and in \cite{L1} coincide. Also, the stratification is normal and balanced in this case.

The construction in \cite{L1} uses heavily representation theoretic tools
like the Frobenius splitting for quantum groups at a root of unity. 
One of the aims of the theory of Seshadri stratification is to provide a purely geometric construction of a standard monomial theory for Schubert varieties, which is also the main point of this article.

\subsection{New geometric approach} 

We finally come to the content of the present paper. Our principal aim is the construction of the standard monomial theory for $X(\tau)$ via the Seshadri stratification by its Schubert subvarieties as recalled above. 

We show that the set $\LS^+_d(\lambda)$ of LS-paths of degree $d$ of shape $\lambda$ is contained in the subset $\Gamma_d$ of values of $\mathcal{V}$ of degree $d$. This follows by a general criterion proved in \cite{CFL} linking our theory to the classical Newton-Okounkov theory. But by the combinatorial character formula for the Demazure module $V(\lambda)_\tau$ proved by the third author in \cite{L2}, the cardinality of $\LS^+_d(\lambda)$ is equal to the dimension of $V(d\lambda)^*$ which is nothing else but the cardinality of $\Gamma_d$. So we have the equality $\LS^+(\lambda) = \Gamma$, also degree by degree. In particular the stratification is normal as follows at once by the definition of $\LS^+(\lambda)$.

In order to deduce that the homogeneous coordinate ring $R$ of $X(\tau)$ has a standard monomial theory, i.e. that it is an LS-algebra, we need two improvements of the theory of Seshadri stratifications developed in the main paper \cite{CFL}.

In \cite{CFL} the quasi-valuation $\mathcal{V}$ is constructed \emph{only} for linearizations $\leq^t$ of the partial order on the poset $A$ that are length preserving: if $\ell(p) < \ell(q)$ then $p \leq^t q$. Here we show that all the constructions in that paper can be carried over for an arbitrary $\leq^t$ (see Lemma~\ref{lemma_quasivaluation_extremal_functions}). 
The definition of a balanced stratification in \cite{CFL} requires: (1) the fan of monoids $\Gamma$ is independent of the linearization $\leq^t$ and (2) there exists a basis of $R$ consisting of common representatives of non-zero leaves for all the quasi-valuations $\mathcal{V}_{\leq^t}$ defined in terms of the linearizations $\leq^t$. In Section \ref{subsection_balanced_stratification} we prove that the condition (2) is a consequence of (1), so the definition of balanced stratification can be weakened.

Once the equality $\Gamma = \LS^+(\lambda)$ has been established, we deduce that $\Gamma$ is independent of the linearization of the Bruhat order used to define $\mathcal{V}$ since $\LS^+(\lambda)$ has clearly this property. Hence the Seshadri stratification is balanced and there exists a vector space basis $x_\pi$, $\pi\in\LS^+(\lambda)$, that represents the non-zero leaves for the quasi-valuation $\mathcal{V}$ defined with respect to an arbitrary linearization of the Bruhat order. We conclude that the algebra $R$ has as set of generators the functions $x_\pi$, $\pi\in\LS^+_1(\lambda)$ and that in the straightening relations for non-standard monomials only standard monomials that are ``greater'' with respect to any linearization $\leq^t$ may appear. In other word $R$ is an LS-algebra and we have constructed a standard monomial theory for the Schubert varieties.

\section{Seshadri stratifications}\label{A_partially_ordered_collection}

In this Section we shortly present Seshadri stratifications for embedded projective varieties; see \cite{CFL} for details. Besides the recollection, there are two new results proved in this section: (a) the length-preserving condition posed on the linearization used in the definition of the quasi-valuation is removed, see Lemma~\ref{lemma_quasivaluation_extremal_functions} and the discussion afterwards; (b) the definition of a balanced stratification is simplified by deducing one requirement from the other, see Section \ref{subsection_balanced_stratification}.

\subsection{Partially ordered set}\label{subsection_poset}

In the whole paper $(A,\leq)$ denotes a partially ordered finite set, \emph{poset} for short. For $p\in A$, we denote by 
${A_p}$ the subset $\{q\in A\mid q\leq p\}$, it is again a poset.

In the sequel we will consider the $\Q$--vector space $\Q^A$ of functions with rational values on $A$. If $p\in A$ then we denote by $e_p\in\Q^A$ the function with value $1$ on $p$ and $0$ on every other $q\in A\setminus\{p\}$; it is clear that the set of functions $\{e_p \mid p\in A\}$, is a vector space basis of $\Q^A$.

The support of an element $\underline{a}\in\Q^A$ is defined as $\supp\underline{a} := \{p\in A\,|\,\underline{a}(p) \neq 0\}$.

Given a subset $B$ of $A$ we will naturally consider $\Q^B$ as a vector sub-space of $\Q^A$.

A total order $\leq^t$ on $A$ refining the given partial order $\leq$ of $A$ is called a \emph{linearization} of $\leq$.
 Given a linearization, we also denote by $\leq^t$ the following lexicographic ordering on $\mathbb{Q}^A$: $\underline{a} \leq^t \underline{b}$ if either $\underline{a} = \underline{b}$ or, denoting by $p$ the maximal element with respect to $\leq^t$ such that $\underline{a}(p) \neq \underline{b}(p)$, we have $\underline{a}(p) < \underline{b}(p)$. This is a total order on $\Q^A$ compatible with the vector addition.

\subsection{Seshadri stratifications}


We start with a brief recollection on Seshadri stratifications, as introduced in \cite{CFL}.

Throughout the paper we fix $\mathbb{K}$ to be an algebraically closed field. Let $V$ be a finite dimensional vector space over $\mathbb{K}$. The hypersurface defined as the vanishing set of a homogeneous polynomial function $f\in\mathrm{Sym}(V^*)$ will be denoted by $\mathcal{H}_f:=\{[v]\in\mathbb{P}(V)\mid f(v)=0\}$.
\index{$\mathcal{H}_f$, vanishing set of a homogeneous polynomial function $f$}

Let $X\subseteq \mathbb P(V)$ be an embedded projective variety with graded homogeneous coordinate ring $R:=\mathbb K[\hat{X}]$, where $\hat{X}\subset V$ is the cone over $X$. Let $X_p$, $p\in A$, be a collection of projective subvarieties $X_p$ in $X$, indexed by a finite set ${A}$. The set $A$ inherits naturally a partial order $\leq$ defined by: for $p,q\in A$, $p\leq q$ if and only if $X_p\subseteq X_q$. We assume that there exists a unique maximal element $p_{\max}$ in $A$ and that $X_{p_{\max}}=X$.

For each $p\in A$, we fix a homogeneous function $f_{p}\in\mathrm{Sym}(V^*)$ of degree larger or equal to one.

\begin{definition}[{\cite[Definition 2.1]{CFL}}]\label{Defn:SS}
The collection of subvarieties ${X_p}$ and homogeneous functions ${f_p}$ for $p\in A$ is termed a \emph{Seshadri stratification}, if the following conditions are satisfied:\index{Seshadri stratification}
\begin{enumerate}
\item[(S1)] the projective varieties $X_p$, $p\in A$, are smooth in codimension one; if $q<p$ is a covering relation in $A$, then $X_q\subseteq X_p$ is a codimension one subvariety; 
\item[(S2)] for any $p\in A$ and any $q\not\leq p$, the function $f_q$ vanishes on $X_p$;
\item[(S3)] for $p\in A$, the set-theoretical intersection satisfies 
$$\mathcal{H}_{f_p}\cap X_p=\bigcup_{q\text{ covered by }p} X_q.$$
\end{enumerate}
In a Seshadri stratification, the functions ${f_p}$ are called \emph{extremal functions}.
\end{definition}

Throughout the paper, the following lemma will be used often without mention.

\begin{lemma}[{\cite[Lemma 2.2]{CFL}}]\label{Lem:SS}
If the collection of subvarieties ${X_p}$ of $X$ and homogeneous functions ${f_p}$ for $p\in A$ defines 
a Seshadri stratification, then
\begin{enumerate}
    \item[(i)] the function $f_p$ does not identically vanish on $X_p$,
    \item[(ii)] all maximal chains in $A$ have the same length, which coincides with $\dim X$. In particular, the poset $A$ is graded.
\end{enumerate}
\end{lemma}


\begin{definition}\label{definition_lenght_function}
Let $p\in A$. The \emph{length} ${\ell}(p)$ of $p$ 
is the length of a (hence any) maximal chain joining $p$ with a minimal element in $A$. 
\end{definition}
According to the above lemma, the length is well-defined and satisfies $\ell(p)=\dim X_p$.

\begin{remark}\label{remark_substratification}
For a fixed $p\in A$, the collection of varieties $X_q$ and the extremal functions 
$f_q$ for $q\in A_p$ satisfies the conditions (S1)--(S3), and hence defines a Seshadri stratification for $X_p\hookrightarrow  \mathbb P(V)$.
\end{remark}

\begin{remark}\label{Rmk:Extend}
We will consider the affine cones of the subvarieties in a Seshadri stratification. It is useful to extend the notation one step further. For a minimal element $p\in A$, the affine cone $\hat{X}_p\cong \mathbb{A}^1$. We set $\hat{A}:=A\cup\{p_{-1}\}$
\index{$\hat{A}$, extended poset}
with $\hat{X}_{p_{-1}}:=\{0\}\in V$. Since the variety $\hat{X}_{p_{-1}}$\index{$\hat{X}_{p_{-1}}$} is contained in the affine cone $\hat{X}_p$ for any minimal element $p\in A$, the set $\hat{A}$ inherits a poset structure by requiring $p_{-1}$ to be the unique minimal element.
\end{remark}

\subsection{A Hasse diagram with bonds}\label{Hasse}

We associate an edge-coloured directed graph to a Seshadri stratification of a projective variety $X$ consisting of subvarieties $X_p$ and extremal functions $f_p$ for $p\in A$.

The Hasse diagram ${\mathcal G}_A$ of the poset $A$ is a directed graph on $A$ whose edges are covering relations, pointing to the larger element.

For a covering relation $p>q$ in $A$, $\hat X_q$ is a prime divisor in $ \hat X_p$. According to (S1), the local ring $\mathcal O_{\hat X_p,\hat X_q}$ is a discrete valuation ring. Let ${\nu_{p,q}}:\mathcal O_{\hat X_p,\hat X_q}\setminus\{0\}\to\mathbb{Z}$ be the associated valuation. Let $R_p:=\mathbb{K}[\hat{X}_p]$ denote the homogeneous coordinate ring of $X_p$. For $f\in R_p\setminus\{0\}$, the value $\nu_{p,q}(f)$ is the \textit{vanishing multiplicity} of $f$ in the divisor $\hat X_q$. The integer $b_{p,q}:=\nu_{p,q}(f_p)$ will be called the \textit{bond} between $p$ and $q$. By (S3), we have $b_{p,q}\geq 1$.

The Hasse diagram with bonds\index{Hasse diagram with bonds} is the diagram with edges coloured with the corresponding bonds: 
$q\stackrel{b_{p,q}}{\longrightarrow}p$.

\begin{remark}\label{Rmk:Extend2}
We extend the construction to the poset $\hat{A}$ (Remark \ref{Rmk:Extend}) and the associated extended Hasse diagram $\mathcal{G}_{\hat{A}}$. For a minimal element $p\in A$, the bond $b_{p,p_{-1}}$ is defined to be the vanishing multiplicity of $f_{p}$ at $\hat{X}_{p_{-1}}=\{0\}$, which coincides with the degree of $f_{p}$.
\end{remark}

\subsection{A family of higher rank valuations}\label{Sec:HigherRank}
From now on we fix a Seshadri stratification on $X\subseteq\mathbb{P}(V)$. Let $R_p:=\mathbb{K}[\hat{X}_p]$ denote the homogeneous coordinate ring of $X_p$ and $\mathbb{K} (\hat{X}_p)$ the field of rational functions on $\hat{X}_p$.

Let $N$ be the least common multiple of all bonds appearing in $\mathcal{G}_{\hat{A}}$.

To a fixed maximal chain $\mathfrak{C}:p_{\max}=p_r>p_{r-1}>\ldots>p_1>p_0$ in $A$, we associate a higher rank valuation $\mathcal{V}_{\mathfrak{C}}:\mathbb{K}[\hat{X}]\setminus\{0\}\to\mathbb{Q}^{\mathfrak{C}}$ as follows.

First choose a non-zero rational function $g_r:=g\in\mathbb{K}(\hat{X})$ and denote by $a_r$ its vanishing order in the divisor $\hat{X}_{p_{r-1}}\subset\hat{X}_{p_r}$. We consider the following rational function  
$$h:=\frac{g_r^N}{f_{p_r}^{N\frac{a_r}{b_r}}}\in\mathbb{K}(\hat{X}_{p_{r}}),$$
where $b_r:=b_{p_r,p_{r-1}}$ is the bond between $p_r$ and $p_{r-1}$. By \cite[Lemma 4.1]{CFL}, the restriction of $h$ to $\hat{X}_{p_{r-1}}$ is a well-defined non-zero rational function on $\hat{X}_{p_{r-1}}$. Let $g_{r-1}$ denote this rational function.
This procedure can be iterated by restarting with the non-zero rational function $g_{r-1}$ on $\hat{X}_{p_{r-1}}$. The output is a sequence of rational functions 
$$g_{\mathfrak{C}}:=(g_r,g_{r-1},\ldots,g_1,g_0)$$ 
with $g_k\in\mathbb{K}(\hat{X}_{p_k})\setminus\{0\}$. 

Collecting the vanishing orders together, we define a map 
$$\mathcal{V}_{\mathfrak{C}}:\mathbb{K}[\hat{X}]\setminus\{0\}\to\mathbb{Q}^\mathfrak{C}\subseteq\Q^A,$$
$$g\mapsto \frac{\nu_{r}(g_r)}{b_r}e_{p_r}+\frac{1}{N}\frac{\nu_{r-1}(g_{r-1})}{b_{r-1}}e_{p_{r-1}}+\ldots+\frac{1}{N^r}\frac{\nu_{0}(g_{0})}{b_{0}}e_{p_0},$$
where $\nu_k:=\nu_{p_k,p_{k-1}}$ is the discrete valuation on the local ring $\mathcal{O}_{\hat{X}_{p_k},\hat{X}_{p_{k-1}}}$, extended to the fraction field.

\begin{theorem}[{\cite[Proposition 6.10, Theorem 6.16]{CFL}}]
$\mathcal{V}_\mathfrak{C}$ is a $\mathbb{Q}^{\mathfrak{C}}$--valued valuation on $\mathbb{K}[\hat{X}]$ having at most one-dimensional leaves.
\end{theorem}

\subsection{The lattice generated by the image of the valuation}

Let $L^{\mathfrak C}_{\mathcal V}\subseteq \Q^A$ be the sublattice generated by $\mathcal{V}_{\mathfrak C}(\mathbb{K}[\hat{X}]\setminus\{0\})$.

\begin{proposition}[{\cite[Proposition 6.12 and 6.13]{CFL}}]\label{proposition_lattice_generators}
There exist $r+1$ rational functions $F_r,\ldots,F_0\in\mathbb{K}(\hat{X})\setminus \{0\}$ such that
\[
\mathcal V_\mathfrak{C}(F_j) = (\underbrace{0,\ldots,0}_{r-j},1/b_j,*,\ldots,*).
\]
Moreover, for any $(r+1)$--tuple of such functions $F_r,\ldots,F_0$, we have
\[
L^{\mathfrak C}_{\mathcal V} = \langle \mathcal V_\mathfrak{C}(F_r),\ldots,\mathcal V_\mathfrak{C}(F_0)\rangle_\mathbb{Z}.
\]
\end{proposition}

The rational functions $F_r,\ldots,F_0\in \mathbb K(\hat{X})$ used above are far from being unique. But all possible choices have one common 
feature. Let $B_{\mathfrak C}$ be the inverse of the rational $(r+1)\times (r+1)$ matrix having as columns the valuations $\mathcal V_{\mathfrak C}(F_r),
\ldots \mathcal V_{\mathfrak C}(F_0)$, this is a lower triangular matrix. It is of the following form: 
\[
B_\mathfrak{C} = \left(
\begin{matrix}
b_r \\
* & b_{r-1}\\
\vdots & \ddots & \ddots\\
* & \cdots & * &  b_0\\
\end{matrix}
\right).
\]

\begin{proposition}[{\cite[Proposition 6.14]{CFL}}]\label{proposition_B_matrix}
Let $v\in\mathbb Q^{\mathfrak C}$. Then $v\in L^{\mathfrak C}_{\mathcal V}$ if and only if $B_\mathfrak{C} \cdot v\in \mathbb{Z}^{r+1}$.
Moreover, the entries of $B_\mathfrak{C}$ are integers.
\end{proposition}

\subsection{A higher rank quasi-valuation}

In order to have a global object, independent of the maximal chain, we introduce a quasi-valuation by taking the minimum of the valuations $\mathcal{V}_\mathfrak{C}$, with $\mathfrak{C}$ a maximal chain in $A$. We refer to \cite[Section 3.1]{CFL} for the definition and basic properties of quasi-valuations.

Fix a linearization $\leq^t$ of the given order $\leq$ on $A$ and let $\leq^t$ be the associated lexicographic order on $\Q^A$ (see Section \ref{subsection_poset}). We define a map
\[
\begin{array}{rcl}
\mathcal{V}:\mathbb{K}[\hat{X}]\setminus\{0\} & \longrightarrow & \mathbb{Q}^A\\
g & \longmapsto & \min_{\mathfrak{C}}\mathcal{V}_{\mathfrak{C}}(g),
\end{array}
\]
where  $\mathfrak{C}$ runs over all maximal chains in $A$ and the minimum is taken with respect to $\leq^t$ on $\Q^A$. Being the minimum of a family of valuations, $\mathcal{V}$ is a quasi-valuation (see \cite[Lemma 3.4]{CFL}).

Of course the quasi-valuation $\mathcal{V}$ \emph{does} depend on the chosen linearization $\leq^t$; when we need to stress such dependence we will write $\mathcal{V}_{\leq^t}$. However there is one case where the value of $\mathcal{V}$ is independent of the linearization; this case is key for the whole construction.
\begin{lemma}\label{lemma_quasivaluation_extremal_functions}
For each extremal function $f_q$, $q\in A$, we have $\mathcal{V}(f_q) = e_q$.
\end{lemma}
\begin{proof} By \cite[Example 6.8]{CFL} we know that $\mathcal{V}_{\mathfrak{C}}(f_q) = e_q$ as long as $q\in\mathfrak{C} : p_r > \cdots > p_0$. If $q\not\in\mathfrak{C}$, then let $p_k\in\mathfrak{C}$ be minimal such that $p_k > q$ (this clearly exists since $p_r = p_{\rm{max}} > q$). Since $p_{k-1}$ and $q$ are not comparable with respect to the partial order on $A$, $f_q$ vanishes on $X_{p_{k-1}}$ by (S2) in the definition of a Seshadri stratification.

On the other hand, $f_q$ does not vanish on $X_q$ by (S3), hence it does not vanish on $X_{p_k}$ since $X_q\subset X_{p_k}$. Using (compare with \cite[Example 5.4]{CFL}):
\[
(f_q)_{\mathfrak{C}} = (f_q, f_q^N, f_q^{N^2}, \ldots, f_q^{N^{r-k}},\ldots)
\]
it follows that
\[
\mathcal V_{\mathfrak C}(f_q)
=\frac{\nu_{p_{k},p_{k-1}}(f_q^{N^{r-k}})} {N^{r-k}b_{p_{k},p_{k-1}}} e_{p_k} +\sum_{i<k} a_ie_{p_i}
\]
for rational numbers $a_i\in\Q$, $0\leq i\leq k-1$. Since $\nu_{p_k,p_{k-1}}(f_q)>0$ and $p_k > q$, this implies $\mathcal V_{\mathfrak C}(f_q) >^t e_q$ as elements of $\Q^A$. This finishes the proof that $\mathcal V(f_q) = e_q$.
\end{proof}

Note that the previous Lemma is proved in \cite{CFL} (see \cite[Lemma 8.3]{CFL}) \emph{only} for linearizations $\leq^t$ preserving length, i.e. such that $p < ^t q$ whenever $\ell(p) < \ell(q)$. However, all the other results of the that paper regarding the quasi-valuation $\mathcal{V}$ hold true for an arbitrary linearization $\leq^t$ with the same proofs once the previous Lemma has been established.

\vskip 0.2cm

Let $\Gamma:=\{\mathcal{V}(g)\mid g\in\mathbb{K}[\hat{X}]\setminus\{0\}\}\subseteq\mathbb{Q}^A$ be the image of the quasi-valuation. We call $\Gamma$ the \emph{fan of monoids} of $\mathcal{V}$; this name is justified by (5) in the following theorem. For a fixed maximal chain $\mathfrak{C}$ in $A$, we define the subset $\Gamma_{\mathfrak{C}}:=\{\underline{a}\in\Gamma\mid \mathrm{supp}\,\underline{a}\subseteq\mathfrak{C}\}$ of $\Gamma$.

\begin{theorem}[{\cite[Proposition 8.6, Corollary 9.1, Lemma 9.6]{CFL}}]\label{theorem_quasi_valuation_main}
The following hold:
\begin{enumerate}
\item[(i)] $\mathcal{V}$ has at most one-dimensional leaves.
\item[(ii)] The fan of monoids $\Gamma$ is contained in $\mathbb{Q}_{\geq0}^A$.
\item[(iii)] For $g\in\mathbb{K}[\hat{X}]\setminus\{0\}$, $\mathcal{V}(g) = \mathcal{V}_\mathfrak{C}(g)$ if and only if $\supp\mathcal{V}(g)\subseteq\mathfrak{C}$.
\item[(iv)] The quasi-valuation is additive if and only if the supports of both functions are contained in the same maximal chain: for $g,h\in\mathbb{K}[\hat{X}]\setminus\{0\}$, $\mathcal{V}(gh) = \mathcal{V}(g) + \mathcal{V}(h)$ if and only if there exists a maximal chain $\mathfrak{C}$ such that $\supp\mathcal{V}(g),\supp\mathcal{V}(h)\subseteq\mathfrak{C}$.
\item[(v)] $\Gamma$ is the union of the finitely generated monoids $\Gamma_\mathfrak{C}$, where $\mathfrak{C}$ runs over all maximal chains.
\end{enumerate}
\end{theorem}

\subsection{Torus action}
The same proof as \cite[Lemma~6.15]{CFL} shows that
\begin{proposition}\label{proposition_torus_action}
Let $T$ be a torus acting on $\hat{X}$ such that for each $p\in A$
\begin{itemize}
\item[(i)] $T$ stabilizes the cone $\hat{X}_p$,
\item[(ii)] the extremal function $f_p$ is $T$--homogeneous.
\end{itemize}
Then for each $\underline{a}\in\Gamma$ there exists a $T$--homogeneous function $f_{\underline{a}}$ such that $\mathcal{V}(f_{\underline{a}}) = \underline{a}$. Moreover, denoting by $\wt(f_p)$ the $T$--weight  of $f_p$, the $T$--weight $\wt(f_{\underline{a}})$ of $f_{\underline{a}}$ is given by
\[
\wt(f_{\underline{a}}) = \sum_{p\in A}\underline{a}(p)\wt(f_p).
\]
\end{proposition}

\begin{remark}\label{remark_kstar_action}
The previous proposition applies to the action of $\mathbb{K}^*$ by multiplication in $V$ since the cones are clearly stabilized by this multiplication and the extremal functions are homogeneous by the definition of a Seshadri stratification; this is indeed the context of \cite[Lemma~6.15]{CFL}. In particular, the degree of a homogeneous function $f_{\underline{a}}$ satisfying $\mathcal{V}(f_{\underline{a}}) = \underline{a}$ is
\[
\deg(f_{\underline{a}}) = \sum_{p\in A}\underline{a}(p)\deg(f_p).
\]
\end{remark}

\subsection{Normal stratification}

An element $\underline{a}\in\Gamma$ is \emph{decomposable} if either it is zero or there exist $\underline{a}_1,\,\underline{a}_2\,\in\Gamma$ such that $\underline{a} = \underline{a}_1 + \underline{a}_2$ with $\min\supp\underline{a}_1 \geq \max\supp\underline{a}_2$. If $\underline{a}$ is not decomposable, then it is said to be \emph{indecomposable}.

\begin{proposition}[{\cite[Proposition~15.3]{CFL}}]\label{proposition_decomposition} Each $\underline{a}\in\Gamma_\mathfrak{C}$ has a decomposition $\underline{a} = \underline{a}_1 + \underline{a}_2 + \ldots + \underline{a}_n$ with $\underline{a}_1,\underline{a}_2, \ldots,\underline{a}_n\in\Gamma_\mathfrak{C}$ indecomposable such that $\min\supp\underline{a}_j \geq \max\supp\underline{a}_{j+1}$ for each $j=1, 2,\ldots,n-1$.
\end{proposition}

For a maximal chain $\mathfrak{C}$ in $A$, the monoid $\Gamma_\mathfrak{C}$ is called \emph{saturated} if $\mathbb{Z}\Gamma_\mathfrak{C}\,\cap\,\Q^\mathfrak{C}_{\geq 0} = \Gamma_\mathfrak{C}$. If, for each maximal chain $\mathfrak{C}$, the monoid $\Gamma_\mathfrak{C}$ is saturated then we say that the stratification is \emph{normal}.

\begin{proposition}[{\cite[Proposition~15.4]{CFL}}]\label{proposition_unique_decomposition} If the stratification is normal then each element of $\Gamma$ has a unique decomposition into indecomposables.
\end{proposition}

\subsection{Balanced stratification}\label{subsection_balanced_stratification}

The quasi-valuation $\mathcal{V}$ depends on the choice of a linearization $\leq^t$ of the partial order $\leq$ of $A$. In particular the fan of monoids $\Gamma$ depends on $\leq^t$; to emphatize this dependence we write $\Gamma_{\leq^t}$. Now denote by $\mathcal{F}$ a family of linearizations of $\leq$. We say that a Seshadri stratification is $\mathcal{F}$--\emph{balanced} if $\Gamma_{\leq^t_1} = \Gamma_{\leq^t_2}$ for each pair of linearizations $\leq^t_1,\leq^t_2\in\mathcal{F}$; we will call this common fan of monoids \emph{the} fan of monoids (with respect to $\mathcal{F}$).

We stress that in \cite{CFL} the notion of a balanced stratification is introduced with respect to the family of all length preserving linearizations (see the discussion after Lemma~\ref{lemma_quasivaluation_extremal_functions}). Moreover the definition in \cite{CFL} (see \cite[Definition 15.7]{CFL}) seems stronger since it requires the existence of a common leaf basis for each quasi-valuation defined in terms of a length preserving linearization. 
Note however that such a basis exists always, see Theorem~\ref{theorem_common_leaf_basis} below.

\begin{lemma}\label{ordercomparison}
Let $\leq_1^t$ and $\leq_2^t$ be linearizations of the partial order $\leq$ on $A$.  If $f\in \mathbb{K}[\hat{X}]\setminus\{0\}$, then $\mathcal V_{\leq^t_1}(f) \leq^t_1 \mathcal V_{\leq^t_2}(f)$.  In particular if $\mathcal V_{\leq^t_1}(f)\neq\mathcal V_{\leq^t_2}(f)$ then $\mathcal V_{\leq^t_1}(f) <^t_1 \mathcal V_{\leq^t_2}(f)$.
\end{lemma}
\begin{proof} By the definition of $\mathcal V_{\leq^t_2}(f)$, there exists a maximal chain $\mathfrak{C}$ such that 
$\mathcal V_{\leq^t_2}(f) = \mathcal V_{\mathfrak{C}}(f)$. 

Now $\mathcal V_{\leq^t_1}(f)$ is the minimum with respect to $\leq^t_1$ of the values $\mathcal V_{\mathfrak{D}}(f)$ with 
$\mathfrak{D}$ running over the set of all maximal chains. In particular, $\mathcal V_{\leq^t_1}(f)$ is less than or equal to 
$\mathcal V_{\mathfrak{C}}(f) = \mathcal V_{\leq^t_2}(f)$.
\end{proof}

\begin{theorem}\label{theorem_common_leaf_basis}
Suppose that the Seshadri stratification is $\mathcal{F}$--balanced and let $\Gamma$ be the fan of monoids with respect to $\mathcal{F}$, then for each $\underline{a}\in\Gamma$ there exists a function $f_{\underline{a}}$ such that $\mathcal{V}_{\leq^t}(f_{\underline{a}}) = \underline{a}$ for each $\leq^t\in\mathcal{F}$.
\end{theorem}
\begin{proof} Fix a degree $d$ and denote by $\Gamma_d$ the set of elements of $\Gamma$ of degree $d$. The proof
will use several induction procedures. 

The first induction procedure is by increasing induction on the number $M$ of length preserving linearizations 
$\leq^t_1,\ldots,\leq^t_M$ of $\leq$ in $\mathcal{F}$. Since the number of linearizations of $\leq$ is finite, the inductive procedure will prove the theorem.

If $M=1$, then obviously for every ${\underline{a}}\in\Gamma$ there exists a homogeneous function $f_{\underline{a}}\in\mathbb{K}[\hat{X}]\setminus\{0\}$ such that  $\mathcal V_{\leq^t_1}(f_{\underline{a}}) = {\underline{a}}$. So in the following we assume $M>1$, and we assume that for the linearizations $\leq^t_1,\ldots,\leq^t_{M-1}$ of $\leq$
and every ${\underline{a}}\in\Gamma$ there exists a homogeneous function $f_{\underline{a}}\in\mathbb{K}[\hat{X}]\setminus\{0\}$ such that  $\mathcal V_{\leq^t_j}(f_{\underline{a}}) = {\underline{a}}$ for $j=1,\ldots,M-1$.

A second inductive procedure is by reverse induction on ${\underline{a}}\in \Gamma_d$ with respect to $\leq_1^t$.

For the base step of this second inductive procedure let $\underline{a}\in\Gamma_d$ be maximal with respect to $\leq^t_1$ and let $f_{\underline{a}}\in\mathbb{K}[\hat{X}]$ be such that $\mathcal{V}_{\leq^t_1}(f_{\underline{a}}) = \underline{a}$. Fix $2\leq j\leq M$, we have $\underline{a} = \mathcal{V}_{\leq^t_1}(f_{\underline{a}})\leq^t_1 \underline{a}' := \mathcal{V}_{\leq^t_j}(f_{\underline{a}})$ by Lemma \ref{ordercomparison}. But, being $\underline{a}'\in\Gamma_d$ and $\underline{a}$ maximal with respect to $\leq^t_1$, we get $\underline{a}' = \underline{a}$. This finishes the proof of the base step.

So, by induction, given $\underline{a}\in\Gamma_d$, we can assume that for all ${\underline{a}}'\in \Gamma_d$, ${\underline{a}}<^t_1{\underline{a}}'$, there exists a homogeneous function of degree $d$: $f_{{\underline{a}}'}\in\mathbb{K}[\hat{X}]\setminus\{0\}$ such that  $\mathcal V_{\leq^t_j}(f_{{\underline{a}}'}) = {\underline{a}}'$ for $j=1,\ldots,M$.

Now by induction on $M$ we know that there exists a function $f_{\underline{a}}$, homogeneous of degree $d$,
such that $\mathcal V_{\leq^t_j}(f_{\underline{a}}) = {\underline{a}}$ for $j=1,\ldots,M-1$. If in addition $\mathcal V_{\leq^t_M}(f_{\underline{a}}) = {\underline{a}}$, then we are done. 

Otherwise we have $\mathcal V_{\leq^t_M}(f_{\underline{a}}) = {\underline{a}}'\not={\underline{a}}$. We start a third inductive procedure, this time by increasing induction with respect to $<^t_M$. We construct a new homogeneous function $\tilde f_{\underline{a}}$ with the property $\mathcal V_{\leq^t_j}(\tilde f_{\underline{a}}) = {\underline{a}}$ for $j=1,\ldots,M-1$ and ${\underline{a}}' <^t_M \mathcal V_{\leq^t_M}(\tilde f_{\underline{a}}) \le^t_M {\underline{a}}$.

Since the number of elements in $\Gamma_d$ is finite, after repeating the procedure a finite number of times we find a homogeneous function $\tilde f_{\underline{a}}$ of degree $d$ with the desired property: $\mathcal V_{\leq^t_j}(\tilde f_{\underline{a}}) = {\underline{a}}$ for $j=1,\ldots,M$, which finishes the proof of the theorem.

It remains to describe the construction of the new function. The assumption $\mathcal V_{\leq^t_M}(f_{\underline{a}}) = {\underline{a}}'\not={\underline{a}}$
implies by Lemma~\ref{ordercomparison} :
$$
\mathcal V_{\leq^t_1}(f_{\underline{a}}) ={\underline{a}} <^t_1{\underline{a}}'=\mathcal V_{\leq^t_M}(f_{\underline{a}})\quad\textrm{and}\quad
\mathcal V_{\leq^t_M}(f_{\underline{a}}) ={\underline{a}}' <^t_M{\underline{a}}=\mathcal V_{\leq^t_1}(f_{{\underline{a}}}).
$$ 
By the reverse induction on $\le^t_1$ we know there exists a function $g_{{\underline{a}}'}$ 
such that $\mathcal V_{\leq^t_j}(g_{{\underline{a}}'}) = {\underline{a}}'$ for $j=1,\ldots,M$. Note that 
$f_{\underline{a}}$ and $g_{{\underline{a}}'}$ are linearly independent because 
$\mathcal V_{\leq^t_1}(f_{{\underline{a}}}) = {\underline{a}} <^t_1{\underline{a}}'$. It follows 
for all $\lambda\in\mathbb K^*$: $f_{\underline{a}}-\lambda g_{{\underline{a}}'}\not=0$. 

We know ${\underline{a}} <^t_1{\underline{a}}'$, but the fact  $\mathcal V_{\leq^t_i}(f_{\underline{a}})={\underline{a}}$ for $i=1,\ldots,M-1$
and $\mathcal V_{\leq^t_M}(f_{\underline{a}})={\underline{a}}'$ implies by Lemma~\ref{ordercomparison}: ${\underline{a}}<^t_i {\underline{a}}'$ for 
$i=1,\ldots,M-1$. So by the minimum rule for quasi-valuations we have for all $\lambda\in\mathbb K^*$:
$$
\mathcal V_{\leq^t_j}(f_{\underline{a}}-\lambda g_{{\underline{a}}'})={\underline{a}}.
$$
It remains to check $\mathcal V_{\leq^t_M}(f_{\underline{a}}-\lambda g_{{\underline{a}}'})$. Since the leaves are only one-dimensional,
one can find $\lambda_0\in\mathbb K^*$ such that $\mathcal V_{\leq^t_M}(f_{\underline{a}}-\lambda_0 g_{{\underline{a}}'}) = \tilde {\underline{a}}>^t_M {\underline{a}}'$.

Set $\tilde f_{\underline{a}}= f_{\underline{a}}-\lambda_0 g_{{\underline{a}}'}$, then by the above we have:
\begin{enumerate}
\item[(i)] $\mathcal V_{\leq^t_j}(\tilde f_{\underline{a}})={\underline{a}}$ for $i=1,\ldots,M-1$;
\item[(ii)] $\mathcal V_{\leq^t_M}(\tilde f_{\underline{a}})=\tilde {\underline{a}}$, where  $\tilde {\underline{a}}>^t_M {\underline{a}}'$.
\end{enumerate}
Lemma~\ref{ordercomparison} implies in addition: $\mathcal V_{\leq^t_M}(\tilde f_{\underline{a}})=\tilde {\underline{a}} \leq^t_M {\underline{a}}= \mathcal V_{\leq^t_1}(\tilde f_{\underline{a}})$,
so we get in addition:
\begin{enumerate}
\item[(iii)] ${\underline{a}}'<^t_M \tilde {\underline{a}}\leq^t_M {\underline{a}}$.
\end{enumerate}
So the new function $\tilde f_{\underline{a}}$ satisfies all the conditions for the inductive procedure, which finishes
the proof.
\end{proof}

\subsection{Standard monomial theory}\label{subsection_standard_monomial_theory}

In this subsection we assume that the Seshadri stratification is normal and $\mathcal{F}$--balanced and we denote by $\Gamma$ the fan of monoids with respect to $\mathcal{F}$.

In what follows we will compare elements of $\Gamma$ with respect to all linearizations in $\mathcal{F}$; so we define a symbol for this relation:
\begin{definition}\label{definition_balanced_order}
Let $\underline{a},\,\underline{b}\in\mathbb{Q}^A$. We write $\underline{a}\trianglelefteq_{\mathcal F}\underline{b}$ if $\underline{a}\leq^t\underline{b}$ for each total order $\leq^t\in\mathcal{F}$.
\end{definition}

Let $\mathbb G\subseteq \Gamma$ be the (possibly infinite) set of indecomposable elements in $\Gamma$. By Proposition \ref{proposition_decomposition}, $\mathbb{G}$ is a generating set of $\Gamma$. For each $\underline{a}\in\mathbb G$, we fix a regular function $x_{\underline{a}}\in\mathbb{K}[\hat{X}]$ satisfying $\mathcal V_{\leq^t}(x_{\underline{a}}) = \underline{a}$ for each $\leq^t\in\mathcal{F}$; these functions exist by Theorem \ref{theorem_common_leaf_basis}. Let $\mathbb{G}_R:=\{x_{\underline{a}}\mid \underline{a}\in\mathbb G\}$.

\begin{definition}\label{definition_standard_monomial}
A monomial $x_{\underline{a}_1}\cdots x_{\underline{a}_n}$ with $\underline{a}_1,\ldots,\underline{a}_n\in\mathbb{G}$ is called \emph{standard} if (up to a suitable reordering of the factors)  $\min\supp\underline{a}_j \geq \max\supp\underline{a}_{j+1}$ for each $j$.
\end{definition}

When writing down a standard monomial $x_{\underline{a}_1}\cdots x_{\underline{a}_n}$, it is understood that for each $j$, $\min\supp\underline{a}_j \geq \max\supp\underline{a}_{j+1}$ holds.

By Proposition \ref{proposition_unique_decomposition}, any element $\underline{a}\in \Gamma$ has a unique decomposition $\underline{a} = \underline{a}_1 + \ldots + \underline{a}_n$ into indecomposable elements and by Theorem \ref{theorem_quasi_valuation_main} each quasi-valuation is additive on standard monomials.

Summarizing we have: 

\begin{proposition}\label{proposition_standard_monomial_basis}
\begin{enumerate}
\item[(i)] The set $\mathbb G_R$ is a generating set for $R$.
\item[(ii)] The set of standard monomials in $\mathbb G_R$ is a vector space basis for $R$. 
\item[(iii)]
If $\underline{a} = \underline{a}_1 + \underline{a}_2 + \ldots + \underline{a}_n$ is the decomposition of $\underline{a}\in\Gamma$ into indecomposables,
then the standard monomial $x_{\underline{a}} := x_{\underline{a}_1}\cdots x_{\underline{a}_n}$ is such that $\mathcal V_{\leq^t}(x_{\underline{a}}) = \underline{a}$ for each $\leq^t\in\mathcal F$.
\item[(iv)]
If a monomial $x_{\underline{a}_1}\cdots x_{\underline{a}_n}$ is not standard, then there exists a straightening relation expressing it as a linear combination of standard monomials
\[
x_{\underline{a}_1}\cdots x_{\underline{a}_n} = \sum_h u_h x_{\underline{a}_{h,1}}\cdots x_{\underline{a}_{h,n_h}},
\]
where $u_h\not=0$ only if 
$\underline{a}_1+\ldots+\underline{a}_n \trianglelefteq_{\mathcal F} \underline{a}_{h,1}+\ldots+\underline{a}_{h,n_h}$.
\item[(v)] If in \emph{(iv)} there exists a chain $\mathfrak{C}$ such that $\supp\underline{a}_i\subseteq\mathfrak{C}$ for all $i=1,\ldots,n$, and $\underline{a}'_1 + \cdots + \underline{a}'_m $ is the decomposition of $\underline{a}_1 + \cdots +\underline{a}_n \in \Gamma$ then the standard monomial $x_{\underline{a}'_1}\cdots x_{\underline{a}'_m}$ appears in the right side of the equation in \emph{(iv)} with a non-zero coefficient.
\end{enumerate}
\end{proposition}

\section{A Seshadri stratification for Schubert varieties}\label{section_schubert_varieties}

In this Section we introduce a Seshadri stratification of a Schubert variety by Schubert sub-varieties. For semplicity, we consider only the case of Schubert varieties for a semisimple simply connected algebraic group. The same construction, slightly modified, holds for any symmetrizable Kac-Moody group. For the general case, all proofs and further details can be found in \cite{CFL3}.

\subsection{Schubert varieties} Let $G$ be a semisimple simply connected algebraic group defined over $\mathbb K$. Fix a maximal torus $T\subseteq G$ and let $B$ be a Borel subgroup containing $T$. Let $N(T)$ be the normalizer of $T$ in $G$ and let $W:=N(T)/T$ be the Weyl group. If $Q\supseteq B$ is a parabolic subgroup, then $B$ acts on $G/Q$ and the closure of an orbit is called a \emph{Schubert variety}. These orbits are in bijection with the elements of $W/W_Q$, where $W_Q$ is the Weyl group of $Q$; for $\tau\in W/W_Q$ let $n_\tau\in N(T)$ be any representative for $\tau$, we define $C(\tau) := B\, n_\tau\, Q$ and $X(\tau) := \overline{B\,n_\tau\, Q}$, these are the \emph{Schubert cell} and the \emph{Schubert variety} in $G/Q$, respectively, corresponding to $\tau$. Since $B$ has only finitely many orbits in $G/Q$, any $B$-invariant irreducible closed subvariety is a Schubert variety. 

Let $\alpha_1,\ldots,\alpha_n$ be the simple roots according to the choice of $B$ and denote by $s_1,\ldots,s_n\in W$ the corresponding simple reflections in the Weyl group. Let $\ell$ be the length function on the Weyl group $W$; the value $\ell(w)$ is the length of any reduced decomposition of $w$ as a product of simple reflections. 

We often identify $W/W_Q$ with the subset $W^Q\subseteq W$ of representatives in $W/W_Q$ of minimal length. The Weyl group is naturally endowed with a partial order by viewing the pair consisting of $W$ and the simple reflections as a Coxeter system. This partial order is called the Bruhat order on $W$. We get an induced Bruhat order on $W/W_Q$ via the identification with the subset $W^Q\subseteq W$. This poset has $\ide$ as its unique minimal element. We view the length function $\ell(\cdot)$ as a function on $W/W_Q$ as follows: we define
$\ell(\tau)$ for $\tau\in W/W_Q$ to be $\ell(\hat\tau)$, where $\hat\tau\in W^Q$ is the unique minimal representative of $\tau$. This is the same as the length in the graded poset $W/W_Q$.

The partial order and the length function have the following geometric interpretation: for $\tau\in W/W_Q$, the length $\ell(\tau)$ is the dimension $\dim X(\tau)$ of the corresponding Schubert variety, and if $\kappa\in  W/W_Q$ is a second element, then
$ X(\kappa)\subseteq X(\tau)$ if and only if $\kappa\le \tau$ in the Bruhat order. 

Denote by $\Lambda=\Lambda(T)$ the weight lattice and, according to the choice of $B$, let $\Lambda^+$ be the monoid of dominant weights. We associate to a dominant weight $\lambda \in \Lambda^+$ a module $V(\lambda)$ as follows: we start with  the irreducible highest weight  representation $V_{\mathbb C}(\lambda)$ for the complex version $G_{\mathbb C}$ of $G$. Let $\mathfrak g_{\mathbb C}:=\textrm{Lie\,} G_{\mathbb C}$ be its Lie algebra. After fixing a highest weight vector $v_\lambda\in V_{\mathbb C}(\lambda)$, we get an admissible lattice $V_{\mathbb Z}(\lambda)\subseteq V_{\mathbb C}(\lambda)$ by applying the Kostant integral $\mathbb Z$-form $\mathfrak{U}(\mathfrak g_{\mathbb C})_{\mathbb Z}$ of the enveloping algebra of $\mathfrak g_{\mathbb C}$ to the fixed highest weight vector  $v_\lambda\in V_{\mathbb C}(\lambda)$. The tensor product $V_{\mathbb Z}(\lambda)\otimes_{\mathbb Z}\mathbb K$ with the field $\mathbb K$ gives the desired $G$--module $V(\lambda)$ of highest weight $\lambda$. 

Let $Q\subseteq G$ be the standard parabolic subgroup of $G$ normalizing the line $\mathbb Kv_\lambda$ through the fixed highest weight vector, i.e. $Q$ is generated by $B$ and the root subgroups $U_{-\alpha}$ for all simple roots $\alpha$ such that $\langle\lambda,\alpha^\vee\rangle=0$. The action of $G$ on $V(\lambda)$ induces an embedding $G/Q\hookrightarrow \mathbb P(V(\lambda))$.

 Fix $X_i$ (resp. $X_{-i}$), $1\leq i\leq n$, to be the Chevalley generator of weight $\alpha_i$ (resp. $-\alpha_i$). For $k\geq 0$, the divided power of $X_{\pm i}$ will be denoted by $X_{\pm i}^{(k)}$.  For $\tau\in W/W_Q$ let $\tau = s_{i_1}\cdots s_{i_t}$ be a reduced decomposition. We associate to such a decomposition the extremal weight vector $v_{\tau}=X_{-i_1}^{(m_1)}\cdots X_{-i_t}^{(m_t)}v_\lambda \in V_{\mathbb{Z}}(\lambda)$ of weight $\tau(\lambda)$, where $m_j=\langle s_{i_{j+1}}\cdots s_{i_t}(\lambda),\alpha_{i_j}^\vee\rangle$ for $1\leq j\leq t$. By the Verma identities, $v_{\tau}$ does not depend on the choice of the reduced decomposition of $\tau$.

We denote by $V(\lambda)_\tau\subseteq V(\lambda)$ the $B$--submodule generated by the orbit $B\cdot v_\tau$; this $B$--submodule is called a \emph{Demazure submodule}. The embedding $G/Q\hookrightarrow \mathbb P(V(\lambda))$ gives us an embedding of the Schubert variety $X(\tau)\hookrightarrow  \mathbb P(V(\lambda)_\tau)$, identifying the orbit $B\cdot [v_\tau]\subseteq \mathbb P(V(\lambda))$ with the open and dense Schubert cell $C(\tau)\subseteq X(\tau)$. In particular, the image of $X(\tau)$ in $\mathbb P(V(\lambda))$ is the closure $\overline{B\cdot [v_\tau]}$ of the orbit.

We will use later, usually without mention, the following representation theoretic interpretation of the homogeneous coordinate ring $\mathbb K[\hat{X}(\tau)]$ of $X(\tau)\subseteq \mathbb P(V(\lambda)_\tau)$.

\begin{lemma}[{\cite[Lemma~3.1]{CFL3}}]\label{lemma_coordinate_ring_schubert}
The degree $d$ part $\mathbb{K}[\hat{X}(\tau)]_d$ of the homogeneous coordinate ring $\mathbb{K}[\hat{X}(\tau)]$ is isomorphic to $V(d\lambda)_\tau^*$ as a $B$--module.
\end{lemma}

We will also need the following result about the multiplicities of certain weights in a Demazure module.

\begin{lemma}[{\cite[Corollary~3.4]{CFL3}}]\label{lemma_multiplicity_one}
Let $\sigma\in W^Q$ be covered by $\tau$ with respect to the Bruhat order and let $\beta$ be a positive root such that $\tau = s_\beta\sigma$. Then the dimension of the $T$--weight spaces in $V(\lambda)_\tau$ of weight $\sigma(\lambda) + j\beta$ is $1$ for all $j = 0, 1, \ldots, \langle\sigma(\lambda),\beta^\vee\rangle$.

\end{lemma}

\subsection{The stratification}

As before, fix a dominant weight $\lambda$, let $Q\supseteq B$ be the standard parabolic subgroup of $G$ associated to $\lambda$, fix $\tau\in W/W_Q$ and consider the Schubert variety $X(\tau)\hookrightarrow \mathbb P(V(\lambda)_\tau)$ embedded in the projective space on the Demazure module $V(\lambda)_\tau$.

The set $A_{\tau}:=\{\sigma\in W/W_Q\mid \sigma\le \tau\}$, endowed with the Bruhat order, is a poset. It has $\ide$ as its unique minimal element and $\tau$ as its unique maximal element. We associate the Schubert variety $X(\sigma)$, which is a closed subvariety of $X(\tau)$, to the element $\sigma\in A_{\tau}$. So we have a collection of subvarieties $X(\sigma)$ of $X(\tau)$, indexed by  the partially ordered set $A_\tau$ such that $\kappa\le\sigma$ if an only if $X(\kappa)\subseteq X(\sigma)$. In addition, all the subvarieties are smooth in codimension one by \cite[Corollary~3.3]{CFL3}. The covering relations correspond to codimension one subvarieties since the length function gives the dimension of a Schubert variety and this function coincides with the length in Definition~\ref{definition_lenght_function}. 

To get a Seshadri stratification, we need in addition a collection of homogeneous functions $f_\sigma\in \mathbb K[V(\lambda)_\tau]$. 
We have fixed for all $\sigma\in A_\tau$ a $T$-eigenvector $v_{\sigma}\in V(\lambda)_\tau$ of weight $\sigma(\lambda)$. The corresponding weight space is one dimensional, so the vector is unique up to a scalar multiple. Denote by $f_\sigma\in (V(\lambda)_\tau)^*$ the corresponding dual vector, i.e. $f_\sigma$ is a $T$--eigenvector of weight $-\sigma(\lambda)$, and $f_\sigma(v_{\sigma})=1$.

\begin{proposition}[{\cite[Proposition~3.5]{CFL3}}]\label{proposition_seshadri_straficication_schubert}
The collection of subvarieties $X(\sigma)$ and linear functions $f_\sigma$, $\sigma\in A_\tau$,
defines a Seshadri stratification for $X(\tau)$.
\end{proposition}

\section{Quasi-valuation and LS-paths}

In this section we present a new approach to the construction of a standard monomial theory for Schubert varieties.

We fix a dominant weight $\lambda$ and a Schubert variety $X(\tau)$ embedded in $\mathbb{P}(V(\lambda)_\tau)$ as in the previous section and continue to consider the stratification of $X(\tau)$ by Schubert subvarieties $X(\sigma)$ in $X(\tau)$ indexed by the poset $A_\tau$, with extremal functions $f_\sigma$, $\sigma\in A_\tau$. Recall that by fixing a linearization $\leq^t$ of the Bruhat order on $A_\tau$, the Seshadri stratification on $X(\tau)$ gives a quasi-valuation $\mathcal{V}_{\leq^t}$ defined as the minimum with respect to $\leq^t$ of the valuations $\mathcal{V}_\mathfrak{C}$ along the chains $\mathfrak{C}$ of $A_\tau$.

In order to simplify the notation we avoid adding $\lambda$ and $\tau$ to certain sets and functions we are going to define (as long as there is no possible confusion).

We start by defining the weight and the degree of an arbitrary element of $\Q^{A_\tau}$.
\begin{definition}\label{definition_weight_degree}
We define the \emph{degree} of $\underline{a}\in\Q^{A_\tau}$ as
\[
\deg(\underline{a}) := \sum_{\sigma\in A_\tau}\underline{a}(\sigma) \, \in \, \Q,
\]
and its \emph{weight} as
\[
\wt(\underline{a}) := \sum_{\sigma\in A_\tau}\underline{a}(\sigma)\sigma(\lambda) \, \in \, \Lambda\otimes\Q.
\]
\end{definition}

\begin{remark}\label{remark_torus_eigenfunction}
The torus $T$ of $G$ acts on the cones $\hat{X}(\sigma)$ for each $\sigma\in A_\tau$ and the functions $f_\sigma$ are $T$--eigenfunctions of degree $1$. So, fixing a linearization $\leq^t$, by Proposition~\ref{proposition_torus_action} and Remark~\ref{remark_kstar_action}, for each $\underline{a}$ in the fan of monoids $\Gamma_{\tau,\leq^t}$ there exists a homogeneous $T$--eigenfunction $f_{\underline{a}}\in\mathbb{K}[\hat{X}(\tau)]$ such that $\mathcal{V}_{\leq^t}(f_{\underline{a}}) = \underline{a}$ whose $T$--weight is $-\wt(\underline{a})$ and whose degree is $\deg(\underline{a})$. We will use this observation also without explicit mention in the rest of the paper.
\end{remark}

As a first step in understanding the fan of monoids $\Gamma_{\tau,\leq^t}$ associated to the Seshadri stratification in Section \ref{section_schubert_varieties}, we compute the bonds of the stratification.

\begin{proposition}\label{proposition_bonds}
Let $\sigma > \eta$ in $A_\tau$ be a covering and $\beta$ be a positive root such that $s_\beta\eta = \sigma$ then
$b_{\sigma,\eta} = \langle\eta(\lambda),\beta^\vee\rangle$.
\end{proposition}
\begin{proof}
This follows by \cite[Lemma~3.2]{CFL3}.
\end{proof}

\begin{definition}
Given a maximal chain $\mathfrak C:\ \tau=\tau_r>\ldots>\tau_0=\ide$ in $A_\tau$, we define the \emph{Lakshmibai-Seshadri lattice}, \emph{LS-lattice} for short, associated to $\mathfrak{C}$ as follows
\begin{equation}\label{EqLSLattice}
\LS_{\mathfrak C}:=\left\{\underline{a} = \sum_{h = 0}^r a_h e_{\tau_h} \in \mathbb Q^{\mathfrak C}
\left\vert 
{\scriptsize
\begin{array}{r}
b_{\tau_r,\tau_{r-1}}a_r\in\mathbb Z\\
b_{\tau_{r-1},\tau_{r-2}}(a_r+a_{r-1})\in\mathbb Z\\
\vdots\\
b_{\tau_1,\tau_0}(a_r+a_{r-1}+\ldots+a_1)\in\mathbb Z\\ 
a_0 +a_1+\ldots+a_r \in \mathbb Z\\
\end{array} }
 \right.\right\} \subseteq \mathbb Q^{\mathfrak C}.
\end{equation}
Further we define the \emph{LS-monoid} associated to $\mathfrak{C}$ as $\LS^+_{\mathfrak{C}} := \LS_{\mathfrak{C}}\,\cap\,\mathbb{Q}_{\geq 0}^{\mathfrak{C}}$. Recall that we view the vector spaces $\mathbb Q^{\mathfrak C}$ as sub-spaces of $\Q^{A_\tau}$; so, the \emph{fan of LS-monoids} is defined as
\[
\LS^+ := \bigcup_{\mathfrak{C}}\LS^+_\mathfrak{C}\subseteq \Q^{A_\tau},
\]where the union runs over all maximal chains $\mathfrak{C}$ of $A_\tau$. The elements of $\LS^+$ are called \emph{LS-paths}.
\end{definition}
As proved in \cite{Ch} and in \cite{CFL3}, the fan of LS-monoids is just a different presentation of the LS-paths in the path model; this justifies the name we have given to the elements of $\LS^+$. Given $\mu\in\Lambda\otimes\Q$ let $\pi_\mu:[0,1]_{\Q} \ni t \longmapsto t\mu \in\Lambda\otimes\Q$ be the straight path from $0$ to $\mu$ in $\Lambda\otimes\Q$. Denote by $*$ the concatenation of paths.
\begin{proposition}[{\cite[Section 8]{Ch}}, {\cite[Appendix I]{CFL3}}]\label{proposition_path_model}
The map
\[
\underline{a} = \sum_{h = 0}^r a_h e_{\tau_h} \,\longmapsto\,\pi_{a_r\tau_r(\lambda)} * \pi_{a_{r-1}\tau_{r-1}(\lambda)} * \cdots * \pi_{a_0\tau_0(\lambda)}
\]
is a bijection from the set $\LS^+_d$ of elements of $\LS^+$ of degree $d$ to the set $\mathbb{B}_\tau(d\lambda)$ of LS-paths of shape $d\lambda$ with support in $A_\tau$. 
\end{proposition}

We are finally ready to state and prove the main theorem of the paper.
\begin{theorem}\label{theorem_leaves_ls_paths}
Let $\leq^t$ be a linearization of the Bruhat order $\leq$ on $A_\tau$, let $\mathcal{V}_{\leq^t}$ be the quasi-valuation associated to the Seshadri stratification of $X(\tau)$ defined in terms of $\leq^t$. Then the fan of monoids $\Gamma_{\tau,\leq^t}$ of $\mathcal{V}_{\leq^t}$ coincides with the fan of LS-monoids $\LS^+$. In particular the stratification is normal and balanced with respect to the family of all linearizations.
\end{theorem}
\begin{proof}
Let us fix a maximal chain $\mathfrak{C}\,:\,\tau = \tau_r > \tau_{r-1} > \cdots > \tau_0 = \ide$ and set $b_j = b_{\tau_j,\tau_{j-1}}$ for $j = 0,\ldots,r$ for short. Our first claim is that there exist rational functions $\eta_r,\eta_{r-1},\ldots,\eta_0$ such that
\[
\mathcal{V}_\mathfrak{C}(\eta_j) = \frac{1}{b_j}e_{\tau_j} - \frac{1}{b_j}e_{\tau_{j-1}}
\]
for $j = 0,1,2,\ldots,r$ (where for $j=0$ we set $e_{\tau_{-1}} = 0$).

We start by proving the claim for the top element in the chain $\mathfrak{C}$, i.e. for $j = r$. So let $\sigma = \tau_{r-1}$ and $b = b_r$ for short and let $\beta$ be the positive root such that $\tau = s_\beta \sigma$. By Lemma~\ref{lemma_multiplicity_one} the weight $\mu := \sigma(\lambda) - \beta = \tau(\lambda) + (b - 1)\beta$ has multiplicity $1$ in the Demazure module $V(\lambda)_\tau$. Hence there exists a function $g \in \mathbb{K}[\hat{X}(\tau)]$ of degree $1$ such that
\[
\wt(g) = -\mu = -\left(\frac{1}{b}\tau(\lambda) + (1 - \frac{1}{b})\sigma(\lambda)\right).
\]
Note that $g^b$ has $T$--weight
\[
-b\mu = -(\tau(\lambda) + (b - 1)\sigma(\lambda)) = -(\sigma(b\lambda) - b\beta)
\]
and, again by Lemma~\ref{lemma_multiplicity_one} applied to the weight $b\lambda$, the $T$--weight space of weight $b\mu$ in $V(b\lambda)_\tau$ has dimension $1$. So, by Lemma~\ref{lemma_coordinate_ring_schubert}, there exists a unique function of degree $b$ in $\mathbb{K}[\hat{X}(\tau)]$ of $T$--weight $-b\mu$ up to multiplication by a non-zero scalar. Since $f_\tau f_\sigma^{b - 1}$ has degree $b$ and weight $-b\mu$, there exists $u\in\mathbb{K}^*$ such that
\[
g^{b} = u\cdot f_\tau f_\sigma^{b - 1}.
\]
If we set $\eta_r := g / f_\sigma$ we have
\[
\mathcal{V}_\mathfrak{C}(\eta_r) = \frac{1}{b}\mathcal{V}_\mathfrak{C}(g^b) - \mathcal{V}_\mathfrak{C}(f_\sigma)  =\frac{1}{b}e_\tau + (1-\frac{1}{b})e_\sigma - e_\sigma = \frac{1}{b}e_\tau - \frac{1}{b}e_\sigma
\]
and our claim is proved for the top element in $\mathfrak{C}$.

Now let $0\leq j < r$. By Remark~\ref{remark_substratification}, the collection of Schubert varieties $X(\sigma)$ and extremal functions $f_{\sigma|\hat{X}(\tau_j)}$ with $\sigma \in A_{\tau_j}$ gives a Seshadri stratification for $X(\tau_j)$. So we consider the maximal chain $\mathfrak{C}_j :\tau_j > \tau_{j-1} > \cdots > \tau_0 = \ide$ in $A_{\tau_j}$ and we apply what already proved for the top element in a chain to $\mathfrak{C}_j$; we get a rational function $\widetilde{\eta}$ on $\hat{X}(\tau_j)$ such that $\mathcal{V}_{\mathfrak{C}_j}(\widetilde{\eta}) = (e_{\tau_j} - e_{\tau_{j-1}})/b_j$.

Let $\eta_j$ be any rational function on $\hat{X}(\tau)$ such that $\eta_{j|\hat{X}(\tau_j)} = \widetilde{\eta}$. Note that $\eta_j$ does not vanish identically on $\hat{X}(\tau_j)$ and is defined on a dense subset of any $\hat{X}(\tau_h)$ for $h = j + 1, \ldots, r$ since the same two properties are true on $\hat{X}(\tau_j) \subset \hat{X}(\tau_h)$. This implies that the first $r-j$ entries of $\mathcal{V}_\mathfrak{C}(\eta_j)$ are $0$ and we conclude that $\mathcal{V}_\mathfrak{C}(\eta_j) = (e_{\tau_j}-e_{\tau_{j-1}})/b_j$ as desired. This finishes the proof of our initial claim.

The functions $\eta_r,\eta_{r-1},\ldots,\eta_0$ we have constructed can be used in Proposition~\ref{proposition_B_matrix} and the related matrix $B_\mathfrak{C}$ is
\[
B_\mathfrak{C} = \left(
\begin{array}{ccccc}
1/b_r \\
-1/b_r & 1/b_{r-1} \\
 & -1/b_{r-1} & 1/b_{r-2} \\
  & & \vdots\\
  & & & -1/b_1 & 1/b_0\\
\end{array}
\right)^{-1} = 
\left(
\begin{array}{ccccc}
b_r \\
b_{r-1} & b_{r-1} \\
b_{r-2} & b_{r-2} & b_{r-2} \\
  & & \vdots\\
b_0  & b_0 & \cdots & b_0 & b_0\\
\end{array}
\right).
\]
So, by Proposition~\ref{proposition_lattice_generators}, we get that the image of the valuation $\mathcal{V}_\mathfrak{C}$ on the field of rational functions on $\hat{X}(\tau)$ is the LS-lattice $\LS_\mathfrak{C}$ associated to the chain $\mathfrak{C}$ (note that $b_0 = 1$).

Now let $\leq^t$ be a linearization of the Bruhat order $\leq$ on $A_\tau$ and let $\mathcal{V}_{\leq^t}$ be the associated quasi-valuation on $\mathbb{K}[\hat{X}(\tau)]$. The monoid $\Gamma_\mathfrak{C}$ of values of $\mathcal{V}_{\leq^t}$ with support in $\mathfrak{C}$ is contained in $\Q^{A_\tau}_{\geq 0}$; hence $\Gamma_\mathfrak{C}\subseteq\LS^+_\mathfrak{C}$. So we see that $\Gamma\subseteq\LS^+$.

If we fix a degree $d\geq 0$ we have that $\Gamma_d$ is a subset of $\LS^+_d$. The cardinality of $\Gamma_d$ is the dimension of the Demazure module $V(d\lambda)_\tau$ by Lemma~\ref{lemma_coordinate_ring_schubert}, but the same dimension is also the cardinality of $\LS_d^+$ by Proposition~\ref{proposition_path_model} and the LS-path character formula (Theorem in Section 5.2 in \cite{L2}). We conclude that $\Gamma_d = \LS^+_d$ for any $d$ and the first statement of the theorem is proved.

In particular the stratification is balanced with respect to the family of all linearizations since $\Gamma$ is independent of the chosen linearization $\leq^t$.

Moreover $\Gamma = \LS^+ = (\bigcup_{\mathfrak{C}}\LS_\mathfrak{C})\cap \Q^{A_\tau}_{\ge 0}$, it is hence clear that $\Gamma$ is saturated and the Seshadri stratification is normal.
\end{proof}

\begin{coro}\label{corollary_geometric_ls_path}
For any LS-path $\underline{a}\in\LS^+$ there exists a homogeneous $T$--eigenfunction $x_{\underline{a}}$ of degree $\deg(\underline{a})$ and weight $-\wt(\underline{a})$ such that
\[
\mathcal{V}_{\leq^t}(x_{\underline{a}}) = \underline{a}
\]
for any linearization $\leq^t$ of the Bruhat order $\leq$ on $A_\tau$.

Hence the combinatorially defined LS-paths are the values of a quasi-valuation, so they encode the vanishing data of regular functions of $\hat{X}(\tau)$ on the net of cones over the Schubert sub-varieties in $X(\tau)$.
\end{coro}

We see a first combinatorial consequence of Theorem \ref{theorem_leaves_ls_paths}. This has already been proved in \cite{De}.
\begin{coro}\label{corollary_gcd_condition}
Let $\sigma = \sigma_d > \sigma_{d-1} > \cdots > \sigma_1 = \eta$ be a maximal chain between $\sigma$ and $\eta$ in $A_\tau$. Then
\[
\gcd(b_{\sigma_d, \sigma_{d-1}}, \ldots, b_{\sigma_2, \sigma_1})
\]
is independent of the chain and depends only on $\sigma$ and $\eta$.
\end{coro}
\begin{proof} Let $\mathfrak{C}$ be any maximal chain in $A_\tau$ containing the given chain. Given a positive integer $n$, by the definition of the LS-lattice, the element
\[
\underline{a} = \frac{1}{n}e_\sigma + \left(1 - \frac{1}{n}\right)e_\eta
\]
is in $\Gamma_{\mathfrak{C}}$ if and only if $n$ divides $\gcd(b_{\sigma_d, \sigma_{d-1}}, \ldots, b_{\sigma_2, \sigma_1})$. If this is the case, then there exists a function $g\in\mathbb{K}[\hat{X}(\tau)]$ such that $\mathcal{V}(g) = \mathcal{V}_\mathfrak{C}(g) = \underline{a}$. But we have $\mathcal{V}_\mathfrak{D}(g) = \underline{a}$ as well for any maximal chain $\mathfrak{D}$ in $A_\tau$ containing $\sigma$ and $\eta$ by (3) of Theorem~\ref{theorem_quasi_valuation_main}. Hence the greatest common divisor does not depend of the chain from $\sigma$ to $\eta$.
\end{proof}

\section{Standard monomial theory}\label{section_SMT}

In the previous section we have determined the fan of monoids $\Gamma_\tau$ for the quasi-valuations of the Seshadri stratification of the Schubert varietiy $X(\tau)\subseteq\mathbb{P}(V(\lambda)_\tau)$ via Schubert sub-varieties. In particular we have obtained a geometric interpretation of the LS-paths as vanishing data of regular functions on the net of Schubert sub-varieties. Moreover we have proved that the Seshadri stratification is normal and balanced with respect to the family $\mathcal{F}$ of all linearizations of the Bruhat order.

We can apply the constructions of the standard monomial theory reviewed in Section \ref{subsection_standard_monomial_theory}. First of all, the indecomposable elements are the LS-paths of degree $1$; this is proved in \cite{Ch}, see also \cite{CFL4}. So, as generating set of the ring $\mathbb{K}[\hat{X}(\tau)]$ we can take the set of functions $x_{\underline{a}}$, with $\underline{a}$ an LS-path of degree $1$, whose existence is stated in Corollary \ref{corollary_geometric_ls_path}. Recall that a monomial $x_{\underline{a}_1}\cdots x_{\underline{a}_n}$ (of degree $n$), with $\underline{a}_1,\ldots,\underline{a}_n\in\LS^+_1$, is \emph{standard} if (up to a suitable reordering of the factors) $\min\supp\underline{a}_j \geq \max\supp\underline{a}_{j+1}$ for each $j$.

\begin{theorem}\label{theorem_standard_monomial_theory_schubert}
\begin{itemize}
    \item[(i)] The set of standard monomials are in bijection with the set $\LS^+$ of LS-paths.
    \item[(ii)] The set of standard monomials is a basis of $\mathbb{K}[\hat{X}(\tau)]$ as a $\mathbb{K}$--vector space.
    \item[(iii)] If the monomial $x_{\underline{a}_1}x_{\underline{a}_2}$ is not standard, then there exists a straightening relation (of degree $2$)
    \[
    x_{\underline{a}_1}x_{\underline{a}_2} = \sum_h u_h x_{\underline{a}_{h,1}}x_{\underline{a}_{h,2}}
    \]
    where $u_h\neq 0$ only if $\underline{a}_1 + \underline{a}_2\trianglelefteq_{\mathcal F} \underline{a}_{h,1}+\underline{a}_{h,2}$ and $\wt(\underline{a}_1 + \underline{a}_2) = \wt(\underline{a}_{h,1} + \underline{a}_{h,2})$.
    \item[(iv)] If in \emph{(iii)} there exists a chain $\mathfrak{C}$ such that $\supp\underline{a}_1,\supp\underline{a}_2\subseteq\mathfrak{C}$ and $\underline{a}'_1 + \underline{a}'_2 $ is the decomposition of $\underline{a}_1 +\underline{a}_2 \in \Gamma$ into indecomposables then the standard monomial $x_{\underline{a}'_1}x_{\underline{a}'_2}$ appears on the right side of the straightening relation in \emph{(iii)} with a non-zero coefficient.
\end{itemize}
\end{theorem}
\begin{proof} This is just a reformulation of Proposition \ref{proposition_standard_monomial_basis}. In (ii) and (iii) we are only considering non-standard monomials of degree $2$, we are using that each $x_{\underline{a}}$, with $\underline{a}\in\LS^+_1$, has degree $1$ and, moreover, that any equation is $T$--homogeneous.
\end{proof}

Comparing the content of the previous Theorem with the definition of LS-algebra in \cite{Ch} and using the Corollary~\ref{corollary_gcd_condition}, we get

\begin{coro}\label{corollary_ls_algebra}
The homogeneous coordinate ring $\mathbb{K}[\hat{X}(\tau)]$ of the Schubert variety $X(\tau)$ embedded in the projective space $\mathbb{P}(V(\lambda)_\tau)$ is an LS-algebra.
\end{coro}

\begin{coro}\label{corollary_relation_degree_two}
Let $ \mathbb{K}[y_{\underline{a}}\,|\,\underline{a}\in\LS^+_1]$ be the polynomial algebra with indeterminates indexed by the LS-paths of degree $1$. The kernel of the map induced by $y_{\underline{a}}\longmapsto x_{\underline{a}}$ is generated by the straightening relations
    \[
    y_{\underline{a}_1}y_{\underline{a}_2} - \sum_h u_h y_{\underline{a}_{h,1}}y_{\underline{a}_{h,2}}
    \]
where $x_{\underline{a}_1}x_{\underline{a}_2}$ is not standard.
\end{coro}
\begin{proof}
This is proved in \cite{Ch} for LS-algebras, see also \cite{CFL4}.
\end{proof}

One of the main outcomes of having an LS-algebra is the existence of a flat degeneration of $\mathbb{K}[\hat{X}(\tau)]$ to a discrete LS-algebra. Note that we have not proved that the LS-algebra in Corollary \ref{corollary_ls_algebra} is \emph{special} (see \cite{Ch}) but this is irrelevant since the same proof of Theorem~11.1 in \cite{CFL} shows that any discrete LS-algebra in the Schubert case is isomorphic to the special discrete LS-algebra (see conjecture stated in \cite[Remark~1]{Ch}).

From the previous corollaries we can derive various 
geometric consequences; for completeness we summarize here the 
most important ones. For part of the results 
(normality, projective normality, Cohen-Macaulay property etc.) there exist various other proofs, see \cite{S4}, or the proofs using Frobenius 
splitting, see \cite{MR}, \cite{RR} and \cite{R}, or the one 
using standard monomial theory \cite{L1}, \cite{LLM}. See also Conclusion~\ref{conclusion} for connections with the present formulation. 
For the details with respect to the following formulation
see \cite{Ch}, \cite{Ch2} and \cite{CFL2}.

\begin{coro}\label{corollary_geometric_consequences}
The embedding of the Schubert variety $X(\tau)$ in $\mathbb{P}(V(\lambda)_\tau)$ is projectively normal and Cohen-Macaulay. It is the intersection of quadrics and degenerates to the reduced union of normal toric varieties, one for each maximal chain in $A_\tau$. This degeneration is compatible with each Schubert sub-variety $X(\sigma)\subseteq X(\tau)$. The degree of the embedding $X(\tau)\subseteq\mathbb{P}(V(\lambda)_\tau)$ is
\[
\sum_{\mathfrak{C}}\prod_{j=1}^r b_{j,\mathfrak{C}},
\]
where $\mathfrak{C}$ runs over all maximal chains in $A_\tau$ and $b_{j,\mathfrak{C}}$, $j=1,\ldots,r$, are the bonds in $A_\tau$ along the chain $\mathfrak{C}$.
\end{coro}

Other results are discussed in \cite{CFL} and \cite{CFL3}; for example: (1). the compatibility of the standard monomial theory to Schubert sub-varieties of $X(\tau)$; (2). the surjectivity of the restriction maps for line bundles to these sub-varieties; (3). surjectivity of the multiplication maps for line bundles; (4). the defining ideal of the sub-varieties in terms of standard monomials.

\begin{conclusion}\label{conclusion}
To conclude, we would like to comment on the possible choices for generators $x_{\underline{a}}$, $\underline{a}$ an LS-path of degree one, on connections with the results on standard monomial  
theory and applications in \cite{LLM} and \cite{L1}, and on the connection to Frobenius splitting and the Demazure character formula. 

The generators $x_{\underline{a}}$, $\underline{a}$ an LS-path of degree one, chosen at the beginning of Section~\ref{section_SMT}, are far from being unique. Though the leaves are one-dimensional, there are many possible choices for a representative of a leaf because we  have only a filtration. 

Sometimes additional arguments can be used to make a preferred choice. An example is the case of the Pl\"ucker embedding of a Grassmann variety $G_{d,n}$ in $\mathbb P(\Lambda^d\mathbb K^n)$. In this case all weight spaces are one dimensional, so, up to scalar multiples, there is a canonical choice for the generators: the Pl\"ucker coordinates, and we are back in the classical situation studied by Hodge, Seshadri and Musili.

The main point in \cite{CFL3} is to show that path vectors $p_{\underline{a}}$ constructed in \cite{L1} are representatives of the corresponding leaves. So looking back, one can say that the method in \cite{L1} provides a tool to pick a representative for a leaf indexed by an LS-path. 

As a consequence, the similarity of the formulas (for example the straightening laws in \cite{L1} and \cite{LLM} and the ones in Theorem~\ref{theorem_standard_monomial_theory_schubert}) is explained by the fact that the results in \cite{L1} and \cite{LLM} are  just a special case of Theorem~\ref{theorem_standard_monomial_theory_schubert}. 
Indeed, the proofs in \cite{LLM} and \cite{L1} hold only for this special choice of the path vectors $p_{\underline{a}}$ as representatives of the corresponding leaves, whereas Theorem~\ref{theorem_standard_monomial_theory_schubert} and the corresponding corollaries are much more general statements: the formulas hold for any choice of representatives of the leaves!

In this context: it would be interesting to know whether other known bases of $V(\lambda)^*_\tau$ are compatible with the filtration and can be used as representatives
$x_{\underline{a}}$ for the leaves associated to an LS-path of degree 1, and hence serve as starting point for a standard monomial theory in the sense of
Theorem~\ref{theorem_standard_monomial_theory_schubert}. It has been shown in the finite type case in \cite{CL} that the elements of the path basis satisfying certain compatibility conditions with respect to a reduced decomposition of the longest element in the Weyl group belong to (up to multiplication by non-zero scalars) the dual canonical basis, so the dual canonical basis of $V(\lambda)^*_\tau$ is a good candidate. he properties of the Mirkovi\'c-Vilonen (MV) basis and its dual basis proved in \cite{BGL}, Section~5.7, suggest that the elements of the dual basis 
could be compatible with the filtration, and hence serve as 
representatives for the leaves as well.

In the present paper we show that all leaves of the quasi-valuation are LS-paths, and combinatorial LS-path character formula implies equality. There are two reasons why it would be interesting to have a proof for the equality without using the character formula: first of 
all, this would give an algebraic geometric proof of the 
combinatorial path character formula in \cite{L1}, and thus provide a geometric proof of a refined character formula as conjectured in \cite{L5} and proved by Kashiwara in terms of quantum groups and crystals in \cite{Kash}.
Secondly, one could start a standard monomial theory as originally suggested in \cite{LS}. In the present approach, the Frobenius splitting of Schubert varieties is incorporated in the construction.
The approach in \cite{LS} has the advantage that it avoids the question of the normality and Frobenius splitting property, and thus serves better as a possible guiding example for other varieties. 

\end{conclusion}

\end{document}